\documentclass[12pt,a4paper]{article}
\usepackage[latin1]{inputenc}
\usepackage[english]{babel}
\usepackage{amsmath}
\usepackage{amsfonts}
\usepackage{amssymb}
\usepackage{indentfirst}
\usepackage{dsfont}
\newtheorem{theorem}{Theorem}[section]
\newtheorem{lemma}[theorem]{Lemma}

\newtheorem{corollary}[theorem]{Corollary}

\newenvironment{definition}[1][Definition]{\begin{trivlist}
\item[\hskip \labelsep {\bfseries #1}]}{\end{trivlist}}

\author{Paul Rochet}
\title{Semiparametric Efficiency of GMM under Approximate Constraints}
\date{}
\newcommand{\mykeywords}{
GMM; Efficiency Bound; Approximate Constraint.}

\begin{document}
\maketitle

\begin{abstract} Generalized empirical likelihood and generalized method of moments are well spread methods of resolution of inverse problems in econometrics. Each method defines a specific semiparametric model for which it is possible to calculate efficiency bounds. By this approach, we provide a new proof of Chamberlain's result on optimal GMM. We also discuss conditions under which GMM estimators remain efficient with approximate moment constraints.
\end{abstract}

\noindent \textbf{Keywords}: \mykeywords\\

\section{Introduction}

We tackle the problem of recovering an unknown probability measure $\mu$ based on a sample $X_1,..,X_n$ of i.i.d. realizations with distribution $\mu$, where additional information on $\mu$ is available in the form of a set of moments equations
\begin{equation}\label{eq0}   \int \Phi(x) d\mu(x) = 0,  \end{equation}
for some vector valued function $\Phi$. This kind of inverse problems finds many practical applications in econometrics, notably when dealing with instrumental variables, see for instance Donald et al. (2009). In some cases, the function $\Phi$ is not known exactly but is assumed to belong to some parametric family $\left\{ \Phi(\theta,.), \theta \in \Theta \subset \mathbb R^d \right\}$. 
We are then interested in the estimation of the true value $\theta_0$ of the parameter, which is, the zero of $\theta \mapsto \int \Phi(\theta,.) d\mu$. The problem of estimating $\theta_0$ in this context has been widely studied in the literature. Two main methods of estimation have been implemented, namely the generalized method of moments (GMM), introduced in Hansen (1982) 
and the generalized empirical likelihood (GEL), developed in Qin and Lawless (1994) 
for this particular context. \\
Although these two methods aim to estimate the same quantity, we point out that they rely on different descriptions of the statistical model. 
Hence, each method is related to a specific semiparametric model, for which we can calculate the efficiency bound for estimating $\theta_0$, following van der Vaart (1998). By this approach, we exhibit necessary conditions for efficiency of GMM and we recover some known results of Hansen (1982) and Chamberlain (1987) on optimal GMM.\\
In many actual situations, the function $\Phi$ may have a complicated form that can only be evaluated numerically. Simulation-based methods have been implemented to deal with approximate constraints, see for instance Mcfadden (1989) and Carrasco and Florens (2000). In this paper we extend the GMM framework to situations where only an approximation $\Phi_m$ of the true constraint function $\Phi$ is available. We provide conditions under which GMM procedures remain efficient asymptotically when replacing $\Phi$ by its approximation.\\
The article falls into the following parts. After exposing the model in Section \ref{sec2}, we make a brief survey on the main methods of estimation in this model, and provide a new proof on the semiparametric efficiency of GMM in Section \ref{sec3}. In Section \ref{sec4}, we discuss the asymptotic efficiency of the method when dealing with an approximate constraint. Proofs are postponed to the Appendix.

\section{The model}\label{sec2}
Let $\mathcal X$ be an open subset of $\mathbb R^q$, endowed with its Borel field $\mathcal B(\mathcal X)$. We observe an i.i.d. sample $X_1,...,X_n$ with unknown distribution $\mu$. We are interested in the estimation of a parameter $\theta_0 \in \Theta \subset \mathbb R^d$ defined by the moment condition
\begin{equation}\label{eq1}  F(\theta_0,\mu):= \int \Phi(\theta_0,x) d\mu(x) = 0,  \end{equation}
where $ \Phi: \Theta \times \mathcal X \to \mathbb R^k$ ($k \geq d$) is a known map. 
The question of estimating efficiently $\theta_0$ relies on the amount of information available on $\mu$. Here, the information given by the moment condition \eqref{eq1} is used to determine the set $\mathcal M$ of possible values for $\mu$ (the model). The true value $\theta_0$ of the parameter being unknown, the distribution of the observations can be any probability measure $\nu$ for which the map $\theta \mapsto F(\theta,\nu)$ is null for some value of $\theta = \theta(\nu) \in \Theta$. The model is therefore defined as
$$ \mathcal M = \left\{ \nu \in \mathcal P(\mathcal X): \exists  \theta=\theta(\nu) \in \Theta: F(\theta,\nu)=0 \right\},    $$ 
where $\theta(\nu)$ is the parameter of interest. In these settings, we aim at calculating the efficiency bound for estimating $\theta$, following van der Vaart (1998). We make the following assumptions ($\Vert . \Vert$ denotes any norm of an Euclidean space). 
\begin{itemize}
\item \textbf{Assumption 1:} $\Theta$ is a compact subset of $\mathbb R^d$.
\item \textbf{Assumption 2:} The map $F(.,\mu)$ is continuous on $\Theta$ and has a unique zero $\theta_0$. Moreover, $\theta_0$ lies in the interior of $\Theta$.
\item \textbf{Assumption 3:} For all $x \in \mathcal X$, the map $\theta \mapsto \Phi(\theta,x)$ is continuous on $\Theta$ and the map $x \mapsto \sup_{\theta \in \Theta} \Vert \Phi(\theta,x) \Vert$ is bounded by some function $\kappa$, integrable with respect to $\mu$.
\item \textbf{Assumption 4:} For all $x \in \mathcal X$, $\theta \mapsto \Phi(\theta,x)$ is twice continuously differentiable in a neighborhood $\mathcal N$ of $\theta_0$. Moreover $\Vert \partial \Phi(\theta,x)/\partial \theta \Vert $ and $\Vert \partial^2 \Phi(\theta,x)/\partial \theta \partial \theta^t \Vert $ are continuous and bounded by an integrable function in this neighborhood ($\partial \Phi(\theta,.)/\partial \theta$ will be noted $\nabla \Phi(\theta,.)$ in the sequel).
\item \textbf{Assumption 5:} The matrices $ D := \int \nabla \Phi(\theta_0,x) d \mu(x) \in \mathbb R^{d\times k}$ and $V:=\int \Phi(\theta_0,x) \Phi^t(\theta_0,x) d \mu(x) \in \mathbb R^{k \times k}$ are of full rank.
\end{itemize} 
These assumptions are usual conditions for this problem, see for instance Qin and Lawless (1994). They ensure the unicity of the parameter $\theta(\nu)$ (which, we recall, is defined as the zero of $F(.,\nu)$) when $\nu$ is close enough to $\mu$ for the total variation topology, and then allow a proper definition of the parameter of interest in the neighborhood of $\mu$. \\
We can now calculate the efficiency bound for estimating $\theta_0$ in this model. For this, we need the following definitions.
\begin{definition} A model $ \left\{ \mu_t, t \geq 0 \right\}$ with $\mu_0 = \mu$ is \textit{differentiable in quadratic mean} at $\mu$ if there exists $g:\mathcal X \rightarrow \mathbb R$ such that $\int_\mathcal X g^2 \ d \mu < \infty$ and 
$$ \lim_{t\rightarrow 0} \int_\mathcal X \left[ \frac 1 t  \left( \sqrt{\frac{d\mu_t}{d\tau_t}} - \sqrt{\frac{d\mu}{d\tau_t}} \right) - \frac 1 2 g \right]^2 d \tau_t = 0,$$
setting for all $t \geq 0$, $\tau_t = \mu_t + \mu$. 
\end{definition}
The function $g$ is called the \textit{score} of $ \left\{ \mu_t, t \geq 0 \right\}$, it satisfies $\int g d \mu = 0$. In the next definition, for all function $T_n: \mathcal X^n \to \Theta$ of the observations, we denote by $\mathcal L(T_n \vert \ \nu)$ the law of $T_n(X_1,...,X_n)$ assuming that $X_1,...,X_n$ are independent with distribution $\nu$.
\begin{definition} An estimator $\hat \theta = \hat \theta(X_1,...,X_n)$ of a parameter $\theta: \mathcal M \to \Theta$ is \textit{locally Gaussian regular} if for all differentiable submodel $ \left\{ \mu_t, t \geq 0 \right\} \subset \mathcal M$ with $\mu_0 = \mu$ and for all positive sequence $(t_n)_{n \in \mathbb N}$ such that $\sqrt n t_n $ is bounded, 
$ \mathcal L( \sqrt n ( \hat \theta - \theta(\mu_{t_n})) \vert \ \mu_{t_n})$ converges weakly towards a Gaussian distribution as $n \to \infty$.
\end{definition}

In a given model, the efficiency bound for estimating a parameter $\theta_0$ is to be understood as a lower bound for the asymptotic variance of locally Gaussian regular estimators of $\theta_0$. An efficiency bound is calculated by considering Fisher Informations of differentiable submodels. We refer to Bickel et al. (1994) and van der Vaart (1998) for further details. 

\begin{theorem}[Theorem 3, Qin and Lawless (1994)]\label{qin} Suppose that Assumptions 1 to 5 hold. The efficiency bound in this model for estimating $\theta_0$ is
$$   B = \left[ D V^{-1} D^t \right]^{-1}. $$
\end{theorem}
Once we have calculated the efficiency bound in our model, the objective is to build an estimator $\hat \theta$ of $\theta_0$ for which the efficiency bound is reached, at least asymptotically in the sense that
$$ \lim_{n \to \infty} n \ \text{var} (\hat \theta) = B.     $$
In some cases, there may not exist any locally Gaussian regular estimates achieving the bound, see for instance examples in Ritov and Bickel (1990). It may also exist estimators having an asymptotic variance smaller than the efficiency bound, in which case the required regularity conditions are not satisfied, as seen in Chapter 2 in Bickel et al. (1994). Such situations will not occur here, as we assume regularity conditions on the model under which GMM and GEL procedures yield regular estimates.

\section{Estimation of the parameter}\label{sec3}
For this problem, we may adopt two natural, although seemingly different, procedures to estimate $\theta_0$, following Chapter 3 in Bickel et al. (1994). Let $\mu_n = \frac 1 n \sum_{i=1}^n \delta_{X_i}$ denote the empirical distribution, where $\delta$ stands for the Dirac measure.  
\begin{itemize}
\item \textbf{Procedure 1:} Find a "smooth" extension $\overline{\theta}$ of $\theta$ over a larger set $\mathcal P \supseteq \mathcal M$ of probability measures containing the empirical distribution $\mu_n$ and define the estimator as $\hat \theta = \overline{\theta} (\mu_n)$.   
\item \textbf{Procedure 2:} Build an approximation $\hat \mu$ of $\mu$ lying in the model $\mathcal M$ and define the estimator as $\hat \theta = \theta(\hat \mu)$.
\end{itemize}

In the literature, two main methods have been implemented for this problem, each one providing a good illustration of each procedure. 

\subsection{Generalized method of moments}\label{secGMM}

The generalized method of moments (GMM) was introduced in Hansen (1982). 
The method consists in replacing in the moment constraint the true measure $\mu$ by its empirical approximation $\mu_n$. Then, find the value of $\theta$ for which $ F(\theta,\mu_n) = \frac{1}{n} \sum_{i=1}^n \Phi(\theta,X_i)$ is as close as possible to $0$ according to a given euclidean norm of $\mathbb R^k$. Precisely, define for $M$ a symmetric positive definite $k \times k$ matrix and $a \in \mathbb R^k$, $\Vert a \Vert_M^2 = a^t M a$, the GMM estimator $\hat \theta$ of $\theta_0$ associated to the norm $\Vert . \Vert_M$ is given by
$$ \hat \theta  = \underset{\theta \in \Theta}{\text{argmin}} \  \Vert F(\theta,\mu_n) \Vert_M.  $$
In practice, the matrix $M$ may have a dependency in $n$, in which case it is chosen to converge towards a symmetric positive definite matrix. However, replacing the matrix by its limit leads to the same first order asymptotic properties of the estimate, under regularity conditions, as pointed out in Newey and Smith (2004). Here, we will assume for simplicity that $M$ is fixed, this being sufficient for our purposes. \\

The generalized method of moments is a good illustration of the first procedure, as the GMM estimator $\hat \theta$ can be seen as the image of the empirical distribution $\mu_n$ by the function 
$$\overline \theta_M(\nu) = \underset{\theta \in \Theta}{\text{argmin}} \ \Vert F(\theta,\nu) \Vert_M, \ \nu \in \mathcal P,  $$
where $\mathcal P$ is an extension of the original model $\mathcal M$, containing $\mu_n$. For sake of generality, $\mathcal P$ is taken as the set of all probability measures $\nu$ for which $ F(.,\nu) $ can take finite values on $\Theta$. Because $\Theta$ is compact, $\mathcal P$ does not depend on the scaling matrix $M$.\\ 

This procedure may seem inefficient at first. Indeed, extending the parameter over a larger model $\mathcal P$ implicitly increases the size of the model, and thus decreases the information available. To be able to provide an efficient estimation, the extension $\overline \theta_M$ must be "smooth" enough so that differentiable submodels in $\mathcal P$ carry at least as much information as the original model. Basically, we want the efficiency bound $\overline B_M$ for estimating $\overline \theta_M$ over $\mathcal P$ not to be higher than the original bound $B$. Since it obviously can not be lower, the objective is to find an efficient extension, for which $\overline B_M = B$. 

\begin{theorem}\label{GMM} Suppose that Assumption 1 to 4 hold. The efficiency bound for estimating $\overline \theta_M$ in $\mathcal P$ is
$$ \textstyle \overline B_M = \left[ D M D^t \right]^{-1} \left[ D M V M D^t  \right] \left[ D M D^t \right]^{-1}. $$
\end{theorem}
This result was originally shown in Chamberlain (1987), although we propose in the Appendix a different proof, based on modern tools on semiparametric efficiency theory.\\

As expected, the efficiency bound $\overline B_M$ in the extended model $\mathcal P$ is larger than in the original model $\mathcal M$ (see Lemma \ref{vmin} in the Appendix). The asymptotic variance of the GMM estimator is precisely the lower bound $\overline B_M$, as shown in Hansen (1982), which proves the efficiency of the method. The theorem also covers the results of Hansen (1982) and Chamberlain (1987) on optimal GMM for $M = V^{-1}$, leading to an efficiency bound in the extended model that is equal to the original bound $B$ of Theorem \ref{qin}.\\ 

Note that the matrix $V$ is generally unknown, since it depends on both $\mu$ and $\theta_0$. In this case, it is replaced by a consistent estimate $\tilde V$, leading to the same asymptotic properties under regularity conditions. Here again, several approaches are possible.\\ 
In the two-step GMM procedure, the estimate $\tilde V$ is built using a preliminary estimator $\tilde \theta$ of $\theta_0$ obtained by a GMM procedure with known scaling matrix (in general, the identity matrix). As a result, $\tilde \theta$ is not in general asymptotically efficient, however, it is $\sqrt n$-consistent and enables to construct a consistent estimate of $V$.\\ 
Another solution is to minimize simultaneously over $\Theta$
\begin{equation}\label{CUE} \theta \mapsto  F(\theta,\mu_n)^t \hat V^{-1}(\theta) F(\theta,\mu_n),   \end{equation}  
where $\hat V^{-1}(\theta)$ denotes here an arbitrary consistent estimate of $V^{-1}(\theta)$, for all $\theta \in \Theta$. The latter approach was introduced in Hansen et al. (1996) as the continuous updating estimation (CUE).\\

\subsection{Generalized empirical likelihood}
Generalized empirical likelihood (GEL) was first applied to this problem in Qin and Lawless (1994), generalizing an idea of Owen (1991). This method is an application of the first procedure. An estimate $\hat \mu$ of $\mu$ is obtained as an entropic projection (in a general sense defined below) of the empirical measure $\mu_n$ onto the model $\mathcal M$. Hence, the measure $\hat \mu$ is the element of the model that minimizes a given $f$-divergence $\mathcal D_f(\mu_n,.)$ with respect to the empirical distribution. Let us recall some definitions.
\begin{definition} Let $f$ be a strictly convex function with $f(1) = f'(1) = 0$, and let $P,Q$ be two probability measures on $\mathcal X$. The $f$-divergence of $Q$ with respect to $P$ is defined as
$$ \mathcal D_f(P,Q) = \int f \left( \frac{dQ}{dP} \right) dP \ \text { if } Q < \! \! < P, \ \mathcal D_f(P,Q) = \infty \ \text{otherwise}. $$
\end{definition}
A $f$-divergence measures the "closeness" between two probability measures. It is non negative and is null only if $P=Q$. This definition can be extended to sets of measures by noting for $\mathcal S$ a subset of $\mathcal P(\mathcal X)$, 
$$ \mathcal D_f(P,\mathcal S) = \inf_{Q \in \mathcal S}  \mathcal D_f(P,Q).$$
\begin{definition} We call entropic projection of $\nu$ on $\mathcal S$ associated to $f$, an element $\nu^* \in \mathcal S$ such that $ \mathcal D_f(\nu,\mathcal S) = \mathcal D_f(\nu,\nu^*) < \infty$.
\end{definition}
An entropic projection always exists as soon as $\mathcal S$ is closed for the total variation topology and $\mathcal D_f(\nu,\mathcal S)$ is finite. Furthermore, it is unique if $\mathcal S$ is also convex (see Csisz\'ar (1967)).\\

\noindent Setting for a fixed $\theta \in \Theta$, $ \mathcal M_\theta := \left\{ \nu \in \mathcal P(\mathcal X): F(\theta,\nu) = 0 \right\}$, the model can be written as $\mathcal M = \cup_{\theta \in \Theta} \mathcal M_\theta$. Thus, the GEL estimator $\hat \theta = \theta(\hat \mu)$ follows by
$$ \hat \theta = \underset{\theta \in \Theta}{\text{argmin}} \ \mathcal D_f(\mu_n, \mathcal M_\theta). $$
Since $\mathcal M_\theta$ is closed and convex, the entropy $\mathcal D_f(\mu_n, \mathcal M_\theta)$ is reached for a unique measure $\hat \mu(\theta)$ in $\mathcal M_\theta$, provided that $\mathcal D_f(\mu_n, \mathcal M_\theta)$ is finite. Then, it appears that computing the GEL estimator involves a two-step procedure. First, build for each $\theta \in \Theta$, the entropic projection $\hat \mu(\theta)$ of $\mu_n$ onto $\mathcal M_\theta$. Then, minimize $\mathcal D_f(\mu_n, \hat \mu(\theta))$ with respect to $\theta$. Since $\hat \mu(\theta)$ is absolutely continuous w.r.t. $\mu_n$ by construction, minimizing $\mathcal D_f(\mu_n,.)$ reduces to find the proper weights $p_1,...,p_n$ to allocate to the observations $X_1,...,X_n$. This turns into a finite dimensional problem, which can be solved by classical convex optimization tools (see for instance Kitamura (2006)). Finally, the GEL estimator $\hat \theta$ can be expressed as the solution to the saddle point problem
$$ \hat \theta = \underset{\theta \in \Theta}{\text{argmin}} \ \underset{(\lambda_1, \lambda_2) \in \mathbb R^{k+1}}{\text{sup}} \ \lambda_1 - \frac 1 n \sum^n_{i=1} f^*(\lambda_1 + \lambda_2^t \Phi(\theta,X_i)),$$
where $f^*(x) =  \sup_{y} \left\{ xy - f(y) \right\}$ denotes the convex conjugate of $f$. \\

Note that if the choice of the $f$-divergence plays a key role in the construction of the estimator, it has no influence on its asymptotic efficiency. Indeed, Qin and Lawless (1994) show that all GEL estimators are asymptotically efficient, regardless of the $f$-divergence used for their computation. Nevertheless, many situations justify the use of specific $f$-divergences. In its original form, empirical likelihood (EL) estimator in Owen (1991) uses the Kullback entropy $K(.,.)$ as $f$-divergence, pointing out that minimizing $K(\mu_n,.)$ reduces to maximizing likelihood among multinomial distributions. Newey and Smith (2004) remark that a quadratic $f$-divergence leads to the CUE estimator of Hansen et al. (1996). Many choices of $f$-divergence can also be given a Bayesian interpretation, using the maximum entropy on the mean (MEM) approach, as shown in Gamboa and Gassiat (1997). 

\section{Dealing with an approximate constraint}\label{sec4}
In many actual applications, only an approximation of the constraint function is available. This may occur if the moment conditions take complicated forms that can only be evaluated numerically or by simulations. Mcfadden (1989) suggested a method dealing with approximate constraint in a similar situation, introducing the method of simulated moments (see also Carrasco and Florens (2000)). In Loubes and Pelletier (2008) and Loubes and Rochet (2009), the authors study a MEM procedure for linear inverse problems with approximate constraints. Here, we propose to extend the GMM framework to a situation with approximate moment conditions. We assume that we observe a sequence $(\Phi_m(\theta,.))_{m \in \mathbb N}$ of approximate constraints, independent with the original sample $X_1,...,X_n$. We are interested in exhibiting sufficient conditions on the sequence $(\Phi_m(\theta,.))_m$ under which estimating $\theta_0$ with GMM procedures remains efficient when the constraint is replaced by its approximation. We discuss the asymptotic properties of the resulting estimates in a framework where both index $n$ and $m$ simultaneously grow to infinity.\\

\noindent In the sequel, we note $W(\theta)$ the inverse of the covariance matrix of $\Phi(\theta,X)$,
$$ \textstyle W(\theta) = \left[ \int \Phi(\theta,.)\Phi^t(\theta,.) d \mu - \int \Phi(\theta,.) d \mu \int \Phi^t(\theta,.) d \mu \right]^{-1}, \theta \in \Theta,$$
while $\hat W(\theta)$ denotes an arbitrary consistent estimator of $W(\theta)$, built from the observations and the constraint function $\Phi(\theta,.)$. In the same way, $W_m(\theta)$ and $\hat W_m(\theta)$ are defined by replacing $\Phi$ by its approximation $\Phi_m$ in the expressions of $W(\theta)$ and $\hat W(\theta)$ respectively.\\ 
For $E$, an Euclidean space endowed with a norm $\Vert . \Vert$, a function $f: \Theta \to E$ and $\mathcal S   \subseteq \Theta$, note
$$ \Vert f \Vert_\mathcal S = \sup_{\theta \in \mathcal S} \Vert f(\theta) \Vert.  $$
We make the following assumptions, where we recall that $\mathcal N$ is a neighborhood of $\theta_0$ defined in Assumption 4.
 
\begin{itemize}
\item \textbf{Assumption 6:} $\Vert \Phi(.,x) \Vert_\Theta$, $\Vert \nabla \Phi(.,x) \Vert_\mathcal N$ and $\Vert \partial^2 \Phi(.,x)/\partial \theta \partial \theta^t \Vert_\mathcal N$ are dominated by a function $\kappa(x)$ such that $\int \kappa^{12}(x) d \mu(x) < \infty$.
\item \textbf{Assumption 7:} For $(\varphi_m)_{m \in \mathbb N}$ a given sequence tending to infinity, the functions $   \varphi_m  \Vert \Phi_m(.,x) - \Phi(.,x) \Vert_\Theta$ and $ \varphi_m \Vert \nabla \Phi_m(.,x) - \nabla \Phi(.,x) \Vert_\mathcal N$ are dominated by a function $\kappa_m(x)$ such that $\sup_{m}\int \kappa_m^{12}(x) d \mu(x) < \infty$. 
\item \textbf{Assumption 8:} 
The random map $\theta \mapsto \hat W(\theta)$ is differentiable on $\mathcal N$ and $\mathbb E(\sqrt n \Vert \hat W - W \Vert_\Theta)^6$, $\mathbb E (\sqrt n \Vert \nabla \hat W - \nabla W \Vert_\mathcal N)^6$ and $\mathbb E(\varphi_m \Vert \hat W_m - \hat W \Vert_\Theta)^3$ are bounded as $m,n$ range over $\mathbb N$.
\end{itemize}
Approximate GMM estimation consists in minimizing over $\Theta$
$$ \theta \mapsto \hat \xi_m(\theta) = \textstyle \left[ \int \Phi_m^t(\theta,.)d \mu_n \right] \tilde W_m \left[ \int \Phi_m(\theta,.)d \mu_n \right], $$
where $\tilde W_m$ is a random matrix with properties to be specified below. It appears that the accuracy of approximate GMM relies on how close the approximate contrast function $\hat \xi_m$ is to its true value (i.e. when the constraint function is known). In this purpose, the scaling matrix $\tilde W_m$ should be chosen as close as possible to the optimal choice $W_0 = W(\theta_0)$.\\ 

As in the situation where the constraint function is known, the two-step GMM procedure provides a natural way to compute the scaling matrix $\tilde W_m$. First build a preliminary estimator $\tilde \theta_m$, minimizing over $\Theta$ 
$$ \theta \mapsto \tilde \xi_m(\theta) = \textstyle \left[ \int \Phi_m^t(\theta,.)d \mu_n \right] \left[ \int \Phi_m(\theta,.)d \mu_n \right], $$
which corresponds to a GMM procedure with identity scaling matrix. Then, define $\tilde W_m = \hat W_m(\tilde \theta_m)$ which is used as scaling matrix in the contrast function $\hat \xi_m$. The resulting approximate two-step GMM estimator satisfies good asymptotic properties as soon as the approximate function $\Phi_m$ converges fast enough towards $\Phi$, as proved in the following theorem.
\begin{theorem}[Robustness of two-step GMM]\label{rGMM} Denote by $\hat \theta_m$ and $\hat \theta$ the two-step GMM estimators obtained respectively with the constraint functions $\Phi_m$ and $\Phi$. If Assumptions 1 to 8 hold,
$$ n \mathbb E(\Vert \hat \theta_m - \hat \theta \Vert)^2 = O(n \varphi_m^{-2}) + o(1).$$ 
In particular, $\hat \theta_m$ is $\sqrt n$-consistent and asymptotically efficient if $n /\varphi_m^{2}$ tends to zero.
\end{theorem}

In the same way, the CUE procedure can be adapted to the case with approximate constraint. Although, the robustness of CUE with approximate constraint requires slightly stronger assumptions.
\begin{itemize}
\item \textbf{Assumption 9:} $\hat W(.)$ and $W(.)$ are twice continuously differentiable on $\mathcal N$ and $\forall \eta > 0, \mathbb P(\Vert d^2 \hat W/d\theta d\theta^t - d^2 W/d\theta d\theta^t \Vert_\mathcal N> \eta) = o(n^{-1})$. Besides, $\hat W_m(.)$ is differentiable on $\mathcal N$ and $\mathbb E(\varphi_m \Vert \nabla \hat W_m - \nabla \hat W \Vert_\mathcal N)^3 $ is bounded as $m,n$ range over $\mathbb N$.
\end{itemize}
Applying the procedure to the approximate constraint, the approximate CUE estimator follows by minimizing over $\Theta$
$$ \theta \mapsto \hat \zeta_m(\theta) = \textstyle \left[ \int \Phi_m^t(\theta,.)d \mu_n \right] \hat W_m(\theta) \left[ \int \Phi_m(\theta,.)d \mu_n \right].$$
\begin{corollary}[Robustness of CUE]\label{rCUE} Denote by $\hat \theta_m$ and $\hat \theta$ the CUE estimators obtained respectively with the constraint functions $\Phi_m$ and $\Phi$. If Assumptions 1 to 9 hold, 
$$ n \mathbb E(\Vert \hat \theta_m - \hat \theta \Vert)^2 = O(n \varphi_m^{-2}) + o(1).$$ 
In particular, $\hat \theta_m$ is $\sqrt n$-consistent and asymptotically efficient if $n /\varphi_m^{2}$ tends to zero.
\end{corollary}

\section{Appendix}

\subsection{Technical lemmas}
\begin{lemma}\label{vmin} For all symmetric positive-definite matrix $M$, 
$$D M D^t \left[ D M V M D^t \right]^{-1} D M D^t \leq D V^{-1} D^t,$$ 
with equality for $M = V^{-1}$.
\end{lemma}

\noindent \textit{Proof.}  Set $A = V^{1/2} M D^t$, $A [A^t A]^{-1} A^t$ is an orthogonal projection matrix with in particular $A [A^t A]^{-1} A^t \leq I \! d$. The inequality holds after multiplying each term by $D V^{-1/2}$ on the left and $V^{-1/2} D^t$ on the right, proving the result.\\

\begin{lemma}\label{cvu} Let $f: \Theta \to \mathbb R$ be a continuous positive function with a unique zero $\theta_0$ lying in the interior of the compact set $\Theta$ and with positive definite Hessian matrix at $\theta_0$. Assume that $f$ is twice continuously differentiable on a neighborhood $\mathcal N$ of $\theta_0$. Let $(f_n)_{n \in \mathbb N}$ be a sequence of positive random functions, twice continuously differentiable on $\mathcal N$, converging in probability towards $f$. Note $\mathcal H = \partial^2 f/\partial \theta \partial \theta^t$ and $\mathcal H_n = \partial^2 f_n/\partial \theta \partial \theta^t$. Moreover, for all $n \in \mathbb N$, let $(f_{m,n})_{m \in \mathbb N}$ be a sequence of positive random functions converging towards $f_n$ as $m \to \infty$. Denote by $\theta_{m,n}$ and $\theta_n$ a minimizer of $f_{m,n}$ and $f_n$ respectively. If the following conditions are met
\begin{description}
\item i) $ \forall \eta >0$, $\mathbb P(\Vert f_n - f \Vert_\Theta > \eta) = o(n^{-1})$ and $\mathbb P(\Vert \mathcal H_n - \mathcal H \Vert_\mathcal N > \eta) = o(n^{-1})$,
\item ii) the $f_{m,n}$ are differentiable on $\mathcal N$ and $C_1 = \sup_{m,n} \mathbb E (\varphi_m \Vert f_{m,n} - f_n \Vert_\Theta)^p$ and $C_2 = \sup_{m,n} \mathbb E (\varphi_m \Vert \nabla f_{m,n} - \nabla f_n \Vert_\mathcal N)^p$ are finite for a $p>0$ and a sequence $(\varphi_m)_{m \in \mathbb N}$ tending to infinity,
\end{description}
then, there is a constant $K> 0$ such that 
$$ \mathbb E \Vert \theta_{m,n} - \theta_n \Vert^p \leq K \varphi_m^{-p} + o(n^{-1}). $$
\end{lemma}

\noindent \textit{Proof.} By continuity of $\mathcal H$ around $\theta_0$, we may assume without loss of generality that $\mathcal N$ is such that $\mathcal H(\theta)$ has all its eigenvalues larger than some constant $2c>0$ for all $\theta \in \mathcal N$. Note $\rho_n$ the smallest eigenvalue of $\mathcal H_n(\theta)$ as $\theta$ ranges over $\mathcal N$. The uniform convergence of $\mathcal H_n$ on $\mathcal N$ in condition \textit{i)} ensures that $ \mathbb P(\rho_n < c) = o(n^{-1})$. Besides, since $\theta_0$ is the unique zero of $f$ on the compact set $\Theta$, we can find a constant $\eta_1>0$ such that $\theta_n$ lies in $\mathcal N$ as soon as $\Vert f_n - f \Vert_\Theta \leq \eta_1$. Hence, still by condition \textit{i)}, $\mathbb P(\theta_n \notin \mathcal N) = o(n^{-1})$. In the same way, there is a constant $\eta_2 >0$ such that $\mathbb P(\theta_{m,n} \notin \mathcal N) \leq \mathbb P(\Vert f_{m,n} - f \Vert_\Theta > 2\eta_2)$, with
\begin{eqnarray*} \mathbb P(\Vert f_{m,n} - f \Vert_\Theta > 2\eta_2) & \leq & \mathbb P(\Vert f_{m,n} - f_n \Vert_\Theta + \Vert f_n - f \Vert_\Theta > 2\eta_2) \\
& \leq & \mathbb P(\Vert f_{m,n} - f_n \Vert_\Theta  > \eta_2)  + \mathbb P( \Vert f_n - f \Vert_\Theta > \eta_2)\\
& \leq & C_1 (\varphi_m \eta_2)^{-p} + o(n^{-1}),
\end{eqnarray*}
by Chebyshev's inequality. Call $\Omega$ the intersection of the three events $\left\{ \theta_n \in \mathcal N \right\}$, $ \left\{ \theta_{m,n} \in \mathcal N \right\}$ and $\left\{ \rho_n \geq c \right\}$, we get $\mathbb P(\Omega^c) \leq C_1(\varphi_m \eta_2)^{-p} + o(n^{-1})$, where $\Omega^c$ denotes the complementary of $\Omega$. Moreover, on $\Omega$, we have
$$ \Vert \nabla f_{m,n} - \nabla f_n \Vert_\mathcal N \geq \Vert \nabla f_n(\theta_{m,n}) \Vert \geq c \Vert \theta_{m,n} - \theta_n \Vert.   $$
Let $\delta$ be the diameter of $\Theta$, it follows that
\begin{eqnarray*} \mathbb E \Vert \theta_{m,n} - \theta_n \Vert^p & \leq &  c^{-p} \ \mathbb E \Vert \nabla f_{m,n} - \nabla f_n \Vert_\mathcal N^p + \delta^p \ \mathbb P(\Omega^c)\\
& \leq & K \varphi_m^{-p} + o(n^{-1}), 
\end{eqnarray*}
for $K = C_1 \delta^p/\eta_2^p + C_2/c^p$.

\subsection{Proofs}

\noindent \textbf{Proof of Theorem \ref{GMM}:}  
Note $\mathcal T$ the set of bounded functions with zero mean under $\mu$. For any $g \in \mathcal T$ and $t>0$, the measure $\mu_t := (1+tg)\mu$ lies in $\mathcal P$ provided that $t$ is small enough. The path $\left\{ \mu_t, t \geq 0 \right\}$ is thus differentiable with score $g$. \\

\noindent The uniform convergence of $F(.,\mu_t)$ towards $F(.,\mu)$ (which follows from Assumptions 1 and 2) ensures the existence of a minimizer $\theta(t)$ of $F(.,\mu_t)$ continuously close to $\theta_0$ as $t \to 0$ and satisfying the first order condition $\gamma_M(\theta(t), \mu_t) = 0$ where 
$$ \textstyle \gamma_M(\theta,\nu) =  \left[ \int (\nabla \Phi(\theta,.)) d \nu \right] M \left[ \int \Phi(\theta,.) d \nu \right], (\theta,\nu) \in \Theta \times \mathcal P.$$
Under Assumptions 2 to 4, the implicit functions theorem applied to the map $(\theta,t) \mapsto \gamma_M(\theta,\mu_t)$ in a neighborhood of $(\theta_0,0)$ warrants the unicity of the minimum $\theta(t)=\overline \theta_M(\mu_t)$. \\ 

\noindent Note $\dot{l} = (\dot{l}_1,...,\dot{l}_d)^t$ the efficient influence function of $\overline \theta_M$. By a Taylor expansion of $ \Phi_{\theta}$ at $\theta_0$ and using that $\gamma_M(\overline \theta_M(\mu_t),\mu_t) = 0$, we get
$$ \textstyle \left[ \int \nabla \Phi_{\theta_0} d\mu_t \right] M \left[ \left[ \int \Phi_{\theta_0} (1+tg) d\mu \right] + \left[ \int \nabla \Phi_{\theta_0}^t d\mu_t \right](\overline \theta_M(\mu_t) - \theta_0) \right] = o(t).$$ 
Since $\overline \theta_M(\mu_t) - \theta_0 = t \int \dot{l} g d \mu + o(t)$ by definition of $\dot{l}$, we obtain after dividing each term by $t$ and making $t$ tend to zero
$$ \textstyle D M  \left[ \int \Phi_{\theta_0} g d\mu \right] = - D M D^t (\int \dot{l} g d \mu).  $$
Since this holds for all $g \in \mathcal T$, we conclude that
$$ \textstyle \dot{l}(.) = - \left[ D M D^t \right]^{-1} D M \Phi_{\theta_0}(.),  $$
checking beforehand that $\dot{l}$ lies in the closure of $\mathcal T$. The efficiency bound is the variance of $\dot{l}(X)$ which proves the result.\\

\noindent \textbf{Proof of Theorem \ref{rGMM}:} For all $\theta \in \Theta$, let
$$\alpha(\theta)  = \int \Phi(\theta,.) d \mu, \ \beta(\theta) = \int \nabla \Phi(\theta,.) d \mu, \ \gamma(\theta) = \int \frac{\partial^2 \Phi(\theta,.)}{\partial \theta \partial \theta^t} d \mu.$$
Besides, note $\hat \alpha(\theta)$ the empirical estimate of $\alpha(\theta)$ and $\hat \alpha_m(\theta)$ the estimate built with $\Phi_m$ and define $\hat \beta(\theta)$, $\hat \gamma(\theta)$, $\hat \beta_m(\theta)$ and $\hat \gamma_m(\theta)$ analogously. \\
First, prove that $\mathbb E (\Vert \tilde \theta_m - \tilde \theta \Vert^6) = O(\varphi_m^{-6}) + o(n^{-1})$. It suffices to verify the conditions of Lemma \ref{cvu} for $p=6$, taking $f_n = \tilde \xi = \hat \alpha \hat \alpha^t$, $f_{m,n} = \tilde \xi_m = \hat \alpha_m \hat \alpha_m^t$ and $f = \alpha \alpha^t$. In this particular case, we have $ \mathcal H_n = \partial^2 \tilde \xi/\partial \theta \partial \theta^t = 2 \hat \beta \hat \beta^t + 2 \hat \gamma \hat \alpha$ and $ \mathcal H = 2 \beta \beta^t + 2 \gamma \alpha$. \\
\noindent First note that $\mathcal H(\theta_0) = 2 \beta^t(\theta_0) \beta(\theta_0)$ is positive definite by Assumption 5. Furthermore, $\hat \alpha(\theta)$ is asymptotically normal and since $ \Vert \Phi(\theta,.) \Vert $ is dominated by a square integrable function $\kappa $ on $\Theta$, we have, for all $\eta >0$, 
$$ \mathbb P( \Vert \hat \alpha - \alpha \Vert_\Theta \geq \eta)  = o(n^{-1}). $$
By assumption, the same argument holds for $\Vert \hat \beta - \beta \Vert_\mathcal N$ and $\Vert \hat \gamma - \gamma \Vert_\mathcal N$. Condition $i)$ in Lemma \ref{cvu} follows directly, noticing that 
$$ \mathcal H_n - \mathcal H = 2 (\hat \beta - \beta) (\hat \beta + \beta)^t + 2 (\hat \gamma - \gamma) \hat \alpha  + 2 \gamma (\hat \alpha - \alpha).$$
Moreover, $ \Vert \tilde \xi_m - \tilde \xi \Vert_\Theta = \Vert \hat \alpha_m^t \hat \alpha_m - \hat \alpha^t \hat \alpha \Vert_\Theta \leq \Vert \hat \alpha_m + \hat \alpha \Vert_\Theta \Vert \hat \alpha_m - \hat \alpha \Vert_\Theta$, yielding 
$$\mathbb E(\varphi_m \Vert \tilde \xi_m - \tilde \xi \Vert_\Theta)^6 \leq \left[\mathbb E(\varphi_m \Vert \hat \alpha_m - \hat \alpha \Vert_\Theta)^{12}\right]^{\frac 1 2} \left[ \mathbb E(\Vert \hat \alpha_m + \hat \alpha \Vert_\Theta)^{12} \right]^{\frac 1 2} $$
by Cauchy-Schwarz inequality. Thus, $\mathbb E(\varphi_m \Vert \tilde \xi_m - \tilde \xi \Vert_\Theta)^6$ is finite by Assumptions 6 and 7. Since $\nabla \tilde \xi_m = 2 \hat \beta_m \hat \alpha_m$ and $\nabla \tilde \xi = 2 \hat \beta \hat \alpha$, assumptions also warrant that $\mathbb E(\varphi_m \Vert \nabla \tilde \xi_m - \nabla \tilde \xi \Vert_\mathcal N)^6 < \infty$. Lemma \ref{cvu} then gives 
$$\Vert \tilde \theta_m - \tilde \theta \Vert^6 = O(\varphi_m^{-6}) + o(n^{-1}).$$
To show the result, we shall now verify that the conditions of Lemma \ref{cvu} hold for $p=2$ with the functions $f_{m,n} = \hat \xi_m$, $f_n = \hat \xi$, $f = \xi$. We now consider $ \mathcal H_n = 2 \hat \beta \tilde W \hat \beta^t + 2 \hat \gamma \tilde W \hat \alpha$ and $ \mathcal H = 2 \beta W_0 \beta^t + 2 \gamma W_0 \alpha$ where $\tilde W = \hat W(\tilde \theta)$ and $W_0 = W(\theta_0)$. \\

\noindent The Hessian matrix $\mathcal H(\theta_0) = 2 \beta(\theta_0) W_0 \beta^t(\theta_0)$ is positive definite by Assumption 5. For condition $i)$ of Lemma \ref{cvu} to be satisfied, we need that for all $\eta >0$, $\mathbb P(\Vert \tilde W - W_0 \Vert > \eta) = o(n^{-1})$. Since $\mathbb P( \tilde \theta \notin \mathcal N) = o(n^{-1})$, we shall only consider the case where $\tilde \theta \in \mathcal N$. By the triangular inequality, we get $ \Vert \hat W(\tilde \theta) - W_0 \Vert \leq \Vert \hat W(\tilde \theta) - \hat W(\theta_0) \Vert + \Vert \hat W(\theta_0) - W_0 \Vert$ and we use that
$$\mathbb P(\Vert \tilde W - W_0 \Vert > \eta) \leq \mathbb P(\Vert \hat W(\tilde \theta) - \hat W(\theta_0) \Vert > \frac \eta 2) + \mathbb P(\Vert \hat W(\theta_0) - W_0 \Vert > \frac \eta 2). $$
Assumption 8 gives $\mathbb P(\Vert \hat W(\theta_0) - W_0 \Vert > \eta/2) = o(n^{-1})$, using Chebyshev's inequality. Furthermore, $\Vert \hat W(\tilde \theta) - \hat W(\theta_0) \Vert \leq \Vert \nabla \hat W \Vert_\mathcal N \Vert \tilde \theta - \theta_0 \Vert$ for a suitable norm in $\mathbb R^{d \times k \times k}$ and for $K > \mathbb E \Vert \nabla W \Vert_\mathcal N$,
$$ \mathbb P (\Vert \nabla \hat W \Vert_\mathcal N \Vert \tilde \theta - \theta_0 \Vert > \frac \eta 2) \leq \mathbb P ( \Vert \tilde \theta - \theta_0 \Vert > \frac{\eta}{2K}) + \mathbb P(\Vert \nabla \hat W \Vert_\mathcal N  > K) = o(n^{-1}) $$
which ensures condition $i)$ of Lemma \ref{cvu}. Write
$$  \hat \xi_m - \hat \xi  =  (\hat \alpha_m - \hat \alpha)^t \tilde W_m \hat \alpha_m + \hat \alpha_m^t (\tilde W_m - \tilde W) \hat \alpha + (\hat \alpha_m - \hat \alpha)^t \tilde W \hat \alpha  $$
where each term can be controlled using Hölder's inequality, as we have for the middle term
\begin{eqnarray*} \mathbb E(\varphi_m  \Vert \hat \alpha_m^t (\tilde W_m - \tilde W) \hat \alpha \Vert_\Theta)^2 & \leq & \mathbb E(\Vert \hat \alpha_m  \Vert_\Theta \varphi_m \Vert \tilde W_m - \tilde W \Vert \Vert \hat \alpha \Vert_\Theta)^2 \\
& \leq & \left[ \mathbb E(\Vert \hat \alpha_m  \Vert_\Theta \Vert  \hat \alpha \Vert_\Theta)^6 \right]^{\frac 1 3} \left[ \mathbb E(\varphi_m \Vert \tilde W_m - \tilde W \Vert)^3 \right]^{\frac 2 3},
\end{eqnarray*}
for an appropriate norm in $\mathbb R^{k \times k}$ for the matrix $\tilde W_m - \tilde W$. Apply the same procedure for the two other terms, with for instance
\begin{eqnarray*} \mathbb E(\varphi_m  \Vert (\hat \alpha_m - \hat \alpha)^t \tilde W_m \hat \alpha_m \Vert_\Theta)^2 & \leq & \mathbb E( \varphi_m \Vert \hat \alpha_m - \hat \alpha  \Vert_\Theta \Vert \tilde W_m \Vert \Vert \hat \alpha \Vert_\Theta)^2 \\
& \leq & \left[ \mathbb E(\varphi_m  \Vert \hat \alpha_m -\hat \alpha  \Vert_\Theta \Vert  \hat \alpha \Vert_\Theta)^6 \right]^{\frac 1 3} \left[ \mathbb E \Vert \tilde W_m \Vert^3 \right]^{\frac 2 3}.
\end{eqnarray*}
To have $\sup_{m,n} \mathbb E( \varphi_m \Vert \hat \xi_m - \hat \xi \Vert_\Theta)^2< \infty$, it suffices to show $\mathbb E(\varphi_m \Vert \tilde W_m - \tilde W \Vert)^3$ is bounded as $n$ and $m$ range over $\mathbb N$, since the rest follows from the first part of the proof. This is true as soon as $\tilde \theta_m$ and $\tilde \theta$ both lie in $\mathcal N$ as we have on the event $\Omega = \{ \tilde \theta, \tilde \theta_m \in \mathcal N \}$,
\begin{eqnarray*} \Vert \tilde W_m - \tilde W \Vert & \leq & \Vert \hat W_m(\tilde \theta_m) - \hat W(\tilde \theta_m) \Vert + \Vert \hat W(\tilde \theta_m) - \hat W(\tilde \theta) \Vert \\
& \leq & \Vert \hat W_m - \hat W \Vert_\Theta + \Vert \nabla \hat W \Vert_\mathcal N \Vert \tilde \theta_m - \tilde \theta \Vert
\end{eqnarray*}
and the result follows from Assumption 8 and by Cauchy-Scharz inequality, since both $\varphi_m \Vert \tilde \theta_m - \tilde \theta \Vert$ and $\Vert \nabla \hat W \Vert_\mathcal N$ have finite moments of order $6$. Hence,
$$\sup_{n,m \in \mathbb N} \mathbb E (\varphi_m \Vert \hat \xi_m - \hat \xi \Vert_\Theta \mathds 1_\Omega)^{2}< \infty. $$
The same reasoning leads to the same conclusion for $\nabla \hat \xi_m$ on $\mathcal N$, namely
$$\sup_{n,m \in \mathbb N} \mathbb E (\varphi_m \Vert \nabla \hat \xi_m - \nabla \hat \xi \Vert_\mathcal N \mathds 1_\Omega )^{2} < \infty.$$
Following the proof of Lemma \ref{cvu}, we show that the complementary of $\Omega$ occurs with negligible probability as $\mathbb P(\Omega^c) = O(\varphi_m^{-6}) + o(n^{-1})$. Since $\Vert \hat \theta_m - \hat \theta \Vert$ remains bounded on $\Omega^c$, we conclude that $\mathbb E(\Vert \hat \theta_m - \hat \theta \Vert \mathds 1_{\Omega^c})^2 = o(\varphi_m^{-2}) + o(n^{-1})$, yielding
$$ \mathbb E(\Vert \hat \theta_m - \hat \theta \Vert)^2 = O(\varphi_m^{-2}) + o(n^{-1}).$$

\vspace{0.4cm}

\noindent \textbf{Proof of Corollary \ref{rCUE}:} The proof is the same as for Theorem \ref{rGMM}, we show that the conditions of Lemma \ref{cvu} are satisfied for $f_{m,n} = \hat \zeta_m = \hat \alpha_m^t \hat W \hat \alpha_m $, $f_n = \hat \zeta = \hat \alpha^t \hat W \hat \alpha$ and $f = \alpha W \alpha$. Condition $i)$ follows from Assumptions 6 and 9, and $\mathbb E (\varphi_m \Vert \nabla \hat \zeta_m - \nabla \hat \zeta \Vert_\mathcal N)^2$ can be bounded as in the proof of the theorem, using the additional condition that $\mathbb E (\varphi_m \Vert \nabla \hat W_m - \nabla \hat W \Vert_\mathcal N)^3$ is bounded.

\bibliographystyle{plain}
\nocite{*}
\bibliography{GMM2}

\end{document}